\documentclass[12pt]{article}
\usepackage{amsfonts}
\usepackage[plainpages,backref,urlcolor=blue]{hyperref}
\usepackage{tikz}

\title {\bf{On large families of bundles over algebraic surfaces}}

\author{\bf{Cristian Anghel and Nicolae Buruiana}}  
\newtheorem{pb}{Problem}[section]
\newtheorem{conj}[pb]{Conjecture}
\newtheorem{tth}[pb]{Theorem}
\newtheorem{defn}[pb]{Definition}
\newtheorem{rem}[pb]{Remark}
\newtheorem{prop}[pb]{Proposition}
\newtheorem{cor}[pb]{Corollary}

\begin{document}

\maketitle

\bigskip
\noindent
{\small {{\bf ABSTRACT:} The aim of this note is to construct sequences of vector bundles with unbounded rank and discriminant on an arbitrary 
algebraic surface. This problem, on principally polarized abelian varieties with cyclic Neron-Severi group generated by the polarization, was considered by Nakashima  in 
connection with  the Douglas-Reinbacher-Yau conjecture on the Strong Bogomolov Inequality. In particular we show that on any surface, the Strong Bogomolov Inequality $SBI_l$ is false for all $l>4$.} }


\


\tableofcontents

\section{Introduction}

Let $S$ a smooth complex projective algebraic surface. For a vector bundle $E$ on $S$ with Chern classes $c_1$ and $c_2$, the discriminant $\Delta$ is defined by the following formula: 
$$\Delta  =2rc_2-(r-1){c_1}^2 .$$
The well known Bogomolov inequality \cite{bogo} asserts that for semi-stable bundles the discriminant is non-negative. In \cite{na} the author of that paper posed the following problem:

\begin{pb}
 Construct a sequence of $\mu$-stable vector bundles $E_m$ with effectively computable discriminants $\delta_m $ 
such that $\delta_m $ go to infinity with $m$.
\end{pb}
\noindent
in connection with the following conjecture from \cite{dou}:

\begin{conj} \label{conj:dry}
 Let $S$ a simply connected surface with trivial or ample canonical bundle. Then the Chern classes of any stable vector bundle 
with nontrivial moduli space obey the following improved Bogomolov inequality:

\[
 2rc_2-(r-1) {c_1}^2-\frac{r^2}{12}c_2(S)\ge 0 .
\]
 
\end{conj}

The above conjecture is false as proved in \cite{co}, \cite{na} and \cite{na1}.
In \cite{na} and \cite{na1} the counterexamples were constructed using the method of elementary
 transformation on arbitrary surfaces as in the conjecture. In \cite{co} there are two families of 
counterexamples: on $K3$ surfaces using the specific methods for such surfaces from \cite{K} and on 
surfaces in $\mathbb{P}^3$ of degree $d \ge 7$ using Horrocks's theory of monads and methods from \cite{co1}. Also, in \cite{nakatrans} is 
introduced the following 

\begin{defn}
A sequence $E_m$ for $m \geq 1$ of $\mu$-stable vector bundles with the rank $r_m$ and discriminant $\Delta _m$ tending to infinity is called a large family.\\
If $r_m={{\it O}}(m^s)$ and 
$\Delta _m={{\it O}}(m^t)$, the large family is of order 
$(s,t)$. 
\end{defn}
\noindent
and is proved that on abelian varieties with cyclic Neron-Severi group there exists large families 
of orders $(1,2)$ and $(1,3)$.

The aim of this note is to construct large families of vector bundles of various orders on arbitrary 
algebraic surfaces  using the ideas from \cite{A}, \cite{A1} 
and \cite{liq}. In {\bf Section \ref{sec:prelim} } we describe the method of Li-Qin from \cite{liq}. Our main result concerning the construction of large families on arbitrary algebraic surfaces is {\bf Theorem \ref{tth:constr}} and is proved in {\bf Section \ref{sec:constr}}. In particular our method shows that on any surface, the Strong Bogomolov Inequality $SBI_l$ is false for all $l>4$. This fact and the relation with the original Douglas-Reinbacher-Yau Conjecture is described in {\bf Section \ref{sec:bogo}}. The main ingredient for our examples is the Serre construction, the results are valid on 
any surface and the bounds are very explicits. Through the paper we shall denote by $\sim $ the equality up to terms of lower degree and by surface we always mean a smooth projective one over $\mathbb C$.

\section{Notations and preliminaries} \label{sec:prelim}

Let $S$ a smooth projective algebraic surface and $L$ a very ample polarization on $S$.
In \cite{A} for the rank-$2$ case and in \cite{liq} for arbitrary rank, the authors uses the Serre construction 
to obtain $L$-stable vector bundles with prescribed Chern classes, provided that $c_2$ is sufficiently great. 
In fact we have the following existence theorem from \cite{liq} which is a generalization of the similar, rank-$2$ case result, from \cite{A}:

\begin{tth} \label{tth:exist}
Let $L$ eventually rescaled such that $r\cdotp L^2 > K_X \cdotp L $
and \\ $\alpha = (r-1)[1+ max(p_g,h^0(S,{\mathcal{O}}_S(rL-c_1+K_S)))+4(r-1)^2\cdotp L^2]$ \\
$+(r-1)c_1\cdotp L-\frac{r(r-1)}{2}\cdotp L^2$. 
Let $c_2 \ge \alpha$. Then there exists an $L$-stable rank $r$ vector bundle $E$ with Chern classes $c_1$ and $c_2$. \\ 
Moreover we have $h^2(X,ad(E))=0$.
\end{tth}
\noindent
For the convenience of the reader we recall the main steps of the proof. \\
- First of all, for $Z$ a reduced $0$-cycle and two line bundles $L$, $L'$ on $S$, it is well known that there exists a locally free extension in $Ext^1({\mathcal O}_S(L)\otimes {\mathcal I}_Z, {\mathcal O}_S(L') )$ iff $Z$ satisfies the Cayley-Bacharach property with respect to the linear system ${\mathcal O}_S(L-L'+K_S)$. \\
- This fact is generalized by Li-Qin in the following: 
\begin{prop}
Consider $r-1$ line bundles $L_1,...,L_{r-1}$ and $0$-cycles \\ $Z_1,...,Z_{r-1}$ on $S$; 
let $W=\oplus( {\mathcal O}_S(L_i)\otimes {\mathcal I}_{Z_i})$. Then, there is a locally free extension in $Ext^1(W,{\mathcal O}_S(L'))$ iff for any $i=1...(r-1)$, $Z_i$ satisfies the Cayley-Bacharach property with respect to the linear system ${\mathcal O}_S(L_i-L'+K_S)$. 
\end{prop}
- The proposition above is used by choosing $L_i=L$ for all $i's$, $L'=c_1+(r-1)L$ and a convenient $0$-cycle $Z$ of appropriate length which produce a locally free extension $E$ with the desired Chern classes.\\
- The next step is the proof of the stability of $E$ which is a consequence of some properties of the chosen $0$-cycle $Z$.\\
- The fact that $h^2(X,ad(E))=0$, which means that the component of $E$ in the moduli space is smooth, follows again from the properties of the $0$-cycle $Z$ using some ideas from \cite{mo}.\\
\newline 
\noindent
The above {\bf Theorem} \ref{tth:exist} will be used in the next section to produce large families of various orders on arbitrary algebraic surfaces.

\section{The construction of large families} \label{sec:constr}
Let's fix an algebraic surface $S$. In what follows, we shall denote by $L_0$ a fixed polarization on $S$ and by $L=aL_0$ a high multiple of $L_0$. We shall consider the first Chern class of the form $c_1=bL_0$ with $r\cdot a-b$ having a constant value; $r$ will be the rank of the bundle we want to obtain.
\begin{rem}
For $r\cdot a-b$ fixed but higher than a certain bound, the Riemann-Roch formula applied for $\mathcal{O}(rL-c_1+K)$ gives:

$$h^0(S,{\mathcal{O}}_S(rL-c_1+K_S)))=$$ $$ \chi (\mathcal{O}_S)+ \frac{1}{2}(rL-c_1+K_S)(rL-c_1)=$$ $$ \frac{r^2L^2}{2}+\frac{rL(-2c_1+K_S)}{2}+ constant \ terms.$$
\end{rem}\noindent
From the above remark, it follows that for large ranks, 
$$max(p_g,h^0(S,{\mathcal{O}}_S(rL-c_1+K_S)))=$$ $$\frac{r^2L^2}{2}+\frac{rL(-2c_1+K_S)}{2}+ constant \ terms.$$
From the existence {\bf Theorem} \ref{tth:exist} the constant $\alpha $ and hence $c_2$ can be taken as
$$4(r-1)^3a^2{L_0}^2+\frac{r(r-1)a^2{L_0}^2}{2}$$
\noindent
up to terms of lower degree. So, $2rc_2$ and the discriminant $\Delta $ have the following asymptotic formulas:
$$2rc_2 \sim 8r(r-1)^3a^2{L_0}^2+r^2(r-1)a^2{L_0}^2,$$
$$\Delta =2rc_2-(r-1){c_1}^2 \sim 8r(r-1)^3a^2{L_0}^2.$$\noindent
Consider now, the parameters $r$ and $a$ varying asymptotically as:
$$r\sim m^s$$
\noindent
and
$$a \sim m^x$$\noindent
for natural $m$ and real pozitive $s,x$. Therefore,
$$\Delta \sim m^{4s+x}$$
\noindent
and we obtain the following:
\begin{tth} \label{tth:constr}
For any smooth projective algebraic surface and any pair $(s,t)$ of positive real numbers with $t>4s$ there exists a large family of order $(s,t)$, i.e. a sequence of stable vector bundles $E_m$ such that their ranks and discriminants satisfy $r_m={{\it O}}(m^s)$ and 
$\Delta _m={{\it O}}(m^t)$.

\end{tth}

\section{Connection with the Strong Bogomolov Inequality } \label{sec:bogo}

In \cite{nakatrans}, Nakashima introduces the definition of the Strong Bogomolov Inequality $SBI_l$:
\begin{defn}
For a positive real number $l$, a polarized algebraic surface $(S,H)$ satisfy the Strong Bogomolov Inequality $SBI_l$ if there is a positive constant $\sigma$ depending only on the surface and the polarization such that for any stable 
vector bundle $E$ of rank $r$ one has the inequality: $$\Delta(E)\geq r^l \sigma .$$
\end{defn}
Of course, the above inequality is a strengthening of the usual Bogomolov's one $\Delta \geq 0$ and {\bf Conjecture} \ref{conj:dry} 
is  a form of $SBI_2$. \\
\
The following obviously proposition from \cite{nakatrans} explain the relation between the notion of large family and the definition of the Strong Bogomolov Inequality $SBI_l$:
\begin{prop}
Let $S$ an algebraic surface with a large family of order $(s,t)$ such that $t<ls$. Then the Strong Bogomolov Inequality $SBI_l$ is false on $S$.
\end{prop} 
From the proposition above and our {\bf Theorem} \ref{tth:constr} we obtain the following:
\begin{cor}
For any polarized algebraic surface $(S,H)$ and any $l>4$ the Strong Bogomolov Inequality $SBI_l$ is false for $(S,H)$.
\end{cor}
\noindent
{\it Proof}: By {\bf Theorem} \ref{tth:constr} one obtain a large family of order $(s,4s+x)$. For any $l>4$, say $l=4+\epsilon$, we want to fulfill the condition $t<ls$. But this means $x<\epsilon s$ which is obviously achieved for an appropriate $x$. {\bf QED}.

\bigskip

\noindent
Cristian Anghel \\
Department of Mathematics\\
Institute of Mathematics of the Romanian Academy\\
Calea Grivitei nr. 21 Bucuresti Romania\\
email:\textit{Cristian.Anghel@imar.ro}

\

\noindent
Nicolae Buruiana \\
Department of Mathematics\\
Institute of Mathematics of the Romanian Academy\\
Calea Grivitei nr. 21 Bucuresti Romania\\
email:\textit{Nicolae.Buruiana@imar.ro}

\end{document}